\documentclass{gtart}
\usepackage{amsmath,amssymb,diagram}

\newcommand{\defn}[1]{\textit{#1}}
\newcommand{\categ}[1]{\mathsf{#1}}
\newcommand{\trunk}[1]{\mathsf{#1}}
\newcommand{\rmod}[1]{\mathcal{#1}}
\newcommand{\augmod}[1]{\mathcal{I}{#1}}
\newcommand{\iso}{\cong}

\newcommand{\Trunk}[1]{\operatorname{Trunk}(#1)}
\newcommand{\Ext}{\operatorname{Ext}}
\newcommand{\Hom}{\operatorname{Hom}}
\newcommand{\Der}{\operatorname{Der}}
\newcommand{\PDer}{\operatorname{PDer}}

\newcommand{\Id}{\operatorname{Id}}
\newcommand{\As}{\operatorname{As}}

\newcommand{\im}{\operatorname{im}}
\newcommand{\invrack}[1]{#1^*}

\newcommand{\Rho}{P}

\renewcommand{\theenumi}{\roman{enumi}}

\renewcommand{\epsilon}{\varepsilon}

\newtheorem{theorem}{Theorem}[section]
\newtheorem{corollary}[theorem]{Corollary}
\newtheorem{lemma}[theorem]{Lemma}
\newtheorem{proposition}[theorem]{Proposition}

\begin{document}
\title{Rack and quandle homology}
\authors{Nicholas Jackson}
\address{Mathematics Institute, University of Warwick, Coventry CV4 7AL,
United Kingdom}
\email{nicholas@maths.warwick.ac.uk}
\begin{abstract}
A \defn{rack} is a set $X$ equipped with a bijective,
self-right-distributive binary operation, and a quandle is a rack which
satisfies an idempotency condition.

In this paper, we further develop the theory of rack and quandle modules
introduced in \cite{jackson:extrq}, in particular defining a tensor
product $\otimes_X$, the notion of a free $X$--module, and the rack
algebra (or \defn{wring}) $\mathbb{Z}X$.

We then apply this theory to define homology theories for racks and
quandles which generalise and encapsulate those developed by
Fenn, Rourke and Sanderson\cite{frs:trunks, frs:james-bundles};
Carter, Elhamdadi, Jelsovsky, Kamada, Langford and Saito
\cite{cjkls:qcoh-state, ces:twisted};
and Andruskiewitsch, Etingof and Gra\~na \cite{ag:pointed, eg:orc}.
\end{abstract}
\primaryclass{18G60}
\secondaryclass{18G35, 18E10}
\keywords{Racks, quandles, modules, homology, cohomology}
\maketitlepage

\section{Introduction}
A \defn{rack} (or \defn{wrack}) is a set $X$ equipped with an asymmetric
binary operation (often written as exponentiation) such that:
\begin{enumerate}
\item[(R1)] For every $a,b \in X$ there is a unique $c \in X$
  such that $c^b = a$.
\item[(R2)] For every $a,b,c \in X$, the \defn{rack identity} holds:
$$a^{bc} = a^{cb^c}$$
\end{enumerate}
In the first of these axioms, the unique element $c$ may be denoted
$a^{\overline{b}}$, although $\overline{b}$ should not itself be regarded
as an element of the rack.  Association of exponents should be understood
to follow the usual conventions for exponential notation.  In particular,
the expressions $a^{bc}$ and $a^{cb^c}$ should be interpreted as $(a^b)^c$
and $(a^c)^{(b^c)}$ respectively.

A rack which, in addition, satisfies the following idempotency criterion
is said to be a \defn{quandle}.
\begin{enumerate}
\item[(Q)] For every $a \in X$, $a^a = a$.
\end{enumerate}
A rack (or quandle) \defn{homomorphism} is a function $f\co X \to Y$
such that $f(a^b) = f(a)^{f(b)}$ for all $a,b \in X$.  Thus, there exist
categories $\categ{Rack}$ and $\categ{Quandle}$.

For any rack $X$, we may construct the \defn{inverted rack} $\invrack{X}$
which is the set $\{ x^* : x \in X \}$ with rack operation ${x^*}^{y^*} :=
(x^{\overline{y}})^*$.

A detailed exposition may be found in the paper by Fenn and Rourke
\cite{fr:racks-links}.

A \defn{trunk} $\trunk{T}$ is an object analogous to a category, and
consists of a class of \defn{objects} and, for each ordered pair $(A,B)$ of
objects, a set $\Hom_{\trunk{T}}(A,B)$ of \defn{morphisms}.  In
addition, $\trunk{T}$ has a number of \defn{preferred squares}
$$\square[A`B`C`D;f`g`h`k]$$
of morphisms, a concept analogous to that of composition in a category.

Given two trunks $\trunk{S}$ and $\trunk{T}$, a \defn{trunk map} or
\defn{functor} $F\co\trunk{S}\to\trunk{T}$ is a map which assigns to
every object $A$ of $\trunk{S}$ an object $F(A)$ of $\trunk{T}$, and to
every morphism $f\co A \to B$ of $\trunk{S}$ a morphism $F(f)\co F(A) \to
F(B)$ of $\trunk{T}$ such that preferred squares are preserved:
$$\square<1`1`1`1;700`500>[F(A)`F(B)`F(C)`F(D);F(f)`F(g)`F(h)`F(k)]$$
or
$$\square<-1`-1`-1`-1;700`500>[F(A)`F(B)`F(C)`F(D);F(f)`F(g)`F(h)`F(k)]$$
A trunk map of the first kind is said to be \defn{covariant}, while a trunk
map of the second kind is said to be \defn{contravariant}.

For any category $\categ{C}$ there is a well-defined trunk
$\Trunk{\categ{C}}$ which has the same objects and morphisms
as $\categ{C}$, and whose preferred squares are the commutative
diagrams in $\categ{C}$.  In particular, we will consider the case
$\Trunk{\categ{Ab}}$, which we will denote $\trunk{Ab}$ where there
is no ambiguity, and which denotes the trunk of Abelian groups.
Trunks were first introduced and studied by Fenn, Rourke and Sanderson
\cite{frs:trunks}.

Given a rack $X$ we define a trunk $\trunk{T}(X)$ as follows:  Let
$\trunk{T}(X)$ have one object for each element $x \in X$, and for each
ordered pair $(x,y)$ of elements of $X$, a morphism $\alpha_{x,y}\co x \to
x^y$ and a morphism $\beta_{y,x}\co y \to y^x$ such that the squares
$$\square<1`1`1`1;1000`500>[x`x^y`x^z`x^{yz} = x^{zy^z};
  \alpha_{x,y}`\alpha_{x,z}`\alpha_{x^y,z}`\alpha_{x^z,y^z}]
  \quad
  \square<1`1`1`1;1000`500>[y`x^y`y^z`x^{yz} = x^{zy^z};
  \beta_{y,x}`\alpha_{y,z}`\alpha_{x^y,z}`\beta_{y^z,x^z}]$$
are preferred for all $x,y,z \in X$.

A trunk map $A\co\trunk{T}(X)\to\trunk{Ab}$ determines an Abelian group
$A_x$ for each element $x \in X$, and for each ordered pair $(x,y)$ of
elements of $X$, Abelian group homomorphisms $\phi_{x,y}\co A_x \to
A_{x^y}$ and $\psi_{y,x}\co A_y \to A_{x^y}$ such that
\begin{eqnarray*}
\phi_{x^y,z}\phi_{x,y} & = & \phi_{x^z,y^z}\phi_{x,z} \\
\phi_{x^y,z}\psi_{y,x} & = & \psi_{y^z,x^z}\phi_{y,z}
\end{eqnarray*}
for all $x,y,z \in X$.  It may often be convenient to denote such a trunk
map by a triple $(A,\phi,\psi)$.  A \defn{rack module} over $X$ (or an
\defn{$X$--module}) is a (covariant) trunk map $\rmod{A} =
(A,\phi,\psi)\co\trunk{T}(X) \to \trunk{Ab}$ such that
$\phi_{x,y}\co A_x \iso A_{x^y}$ is an isomorphism, and
\begin{equation}
\label{eqn:rmod}
\psi_{z,x^y}(a) = \phi_{x^z,y^z}\psi_{z,x}(a) + \psi_{y^z,x^z}\psi_{z,y}(a)
\end{equation}
for all $a \in A_z$ and $x,y,z \in X$.

If $x,y$ lie in the same orbit of $X$ then this implies that $A_x \iso
A_y$ (although the isomorphism is not necessarily unique).  For racks
with more than one orbit it follows that if $x \not\sim y$ then $A_x$
need not be isomorphic to $A_y$.  Rack modules where the constituent
groups are nevertheless all isomorphic are said to be \defn{homogeneous},
and those where this is not the case are said to be \defn{heterogeneous}.
It is clear that modules over transitive racks (racks with a single orbit)
must be homogeneous.

An $X$--module $\rmod{A}$ of the form $(A,\Id,0)$ (so that $\phi_{x,y} =
\Id\co A_x \to A_{x^y}$ and $\psi_{y,x}$ is the zero map $A_y \to
A_{x^y}$) is said to be \defn{trivial}.

Given two $X$--modules $\rmod{A} = (A,\phi,\psi)$ and $\rmod{B} =
(B,\chi,\omega)$, a \defn{homomorphism} of $X$--modules, or an
\defn{$X$--map}, is a natural transformation $f\co \rmod{A} \to \rmod{B}$.
That is, a collection $f = \{ f_x\co A_x
\to B_x : x \in X \}$ of Abelian group homomorphisms such that
\begin{eqnarray*}
\phi_{x,y}f_x & = & f_{x^y}\phi_{x,y} \\
\psi_{y,x}f_y & = & f_{x^y}\psi_{y,x}
\end{eqnarray*}
for all $x,y \in X$.

As shown in a previous paper \cite{jackson:extrq}, these objects are the
Beck modules in the category $\categ{Rack}$, and hence form an Abelian
category (denoted $\categ{RMod}_X$) suitable for the definition of
homology and cohomology theories.

Specialising to the subcategory $\categ{Quandle}$ yields a related
construct suitable for use in quandle homology and cohomology theories.
A \defn{quandle module} is a rack module $\rmod{A} = (A,\phi,\psi)$
which satisfies the additional criterion
\begin{equation}
\label{eqn:qmod}
\psi_{x,x}(a) + \phi_{x,x}(a) = a
\end{equation}
for all $a \in A_x$ and $x \in X$.  These objects form a category
$\categ{QMod}_X$, which is equivalent to the category of Beck modules in
$\categ{Quandle}$.

This paper contains part of my doctoral thesis \cite{jackson:thesis}.
I am grateful to my supervisor Colin Rourke, and to Alan Robinson, Ronald
Brown, and Simona Paoli for many helpful discussions, comments and advice
over the past few years.

\section{Right $X$--modules}
A contravariant trunk map $\rmod{A} = (A,\phi,\psi) \co \trunk{T}(X)
\to \trunk{Ab}$ determines, as in the covariant case, an Abelian group
$A_x$ for each $x \in X$.  In addition, we obtain homomorphisms
$\phi^{x,y}\co A_{x^y} \to A_x$ and $\psi^{y,x}\co A_{x^y} \to A_y$,
satisfying
$$\phi^{x,y}\phi^{x^y,z} = \phi^{x,z}\phi^{x^z,y^z}$$
and
$$\psi^{y,x}\phi^{x^y,z} = \phi^{y,z}\psi^{y^z,x^z}$$
for all $x,y,z \in X$.  A \defn{right rack module} over $X$ (or a
\defn{right $X$--module}), is such a trunk map in which each $\phi^{x,y}$
is an isomorphism, and
$$\psi^{z,x^y} = \psi^{z,x}\phi^{x^z,y^z} + \psi^{z,y}\psi^{y^z,x^z}$$
for all $x,y,z \in X$.

We may choose to refer to the rack modules earlier
defined\cite{jackson:extrq} as \defn{left} rack modules where necessary.

Given two right $X$--modules $\rmod{A} = (A,\phi,\psi)$ and $\rmod{B} =
(B,\chi,\omega)$, a \defn{homomorphism} (or \defn{$X$--map})
$f\co\rmod{A}\to\rmod{B}$ is a natural transformation.  That is, an
Abelian group homomorphism $f_x\co A_x \to B_x$ for each $x \in X$, such
that
$$f_x\phi^{x,y} = \chi^{x,y}f_{x^y}$$
and
$$f_y\psi^{y,x} = \omega^{y,x}f_{x^y}$$
for all $x,y \in X$.  We denote the category thus formed by $\categ{RMod}^X$.

\begin{proposition}
For any rack $X$, there is a categorical equivalence $\categ{RMod}_X \iso
\categ{RMod}^{X^*}$, where $X^*$ denotes the inverted rack of $X$.
\end{proposition}
\begin{proof}
Let $\rmod{A} = (A,\phi,\psi)$ be an arbitrary left $X$--module.  We may
construct a right $X$--module as follows.  For each $x \in X$, let
$B_{x^*} = A_x$, and define Abelian group homomorphisms
\begin{eqnarray*}
\chi^{x^*,y^*}\co B_{{x^*}^{y^*}} \to B_{x^*} & ; &
  a \mapsto \phi_{x^{\bar{y}},y}(a) \\
\omega^{y^*,x^*}\co B_{{x^*}^{y^*}} \to B_{y^*} & ; &
  a \mapsto \psi_{x^{\bar{y}},y^{y\bar{x}\bar{y}}}(a).
\end{eqnarray*}
A routine calculation confirms that $\rmod{B} = (B,\chi,\omega)$ is a
right $X^*$--module.  Furthermore, given another left $X$--module
$\rmod{C} = (C,\gamma,\eta)$ and a homomorphism $f\co \rmod{A}
\to \rmod{C}$ of left $X$--modules, we may construct another right
$X^*$--module $\rmod{D} = (D,\delta,\theta)$ and a homomorphism $g\co
\rmod{B} \to \rmod{D}$ of right $X^*$--modules such that $g_{x^*} = f_x$
for all $x \in X$.  Another routine calculation verifies that $g$ is a
natural transformation.  This assignment ($\rmod{A} \mapsto \rmod{B}$, $f
\mapsto g$) determines a functor $F\co \categ{RMod}_X \to \categ{RMod}^{X^*}$.

Conversely, given a right $X^*$--module $\rmod{B} = (B,\chi,\omega)$, we
may construct a left $X$--module $\rmod{A} = (A,\phi,\psi)$ by setting
$A_x = B_{x^*}$ and
\begin{eqnarray*}
\phi_{x,y}\co A_x \to A_{x^y} & ; &
  a \mapsto \chi^{{x^*}^{\overline{y^*}},{y^*}}(a) \\
\psi_{y,x}\co A_y \to A_{x^y} & ; & b \mapsto
  \omega^{{x^*}^{\overline{y^*}},{y^*}^{{y^*}\overline{x^*}\,\overline{y^*}}}(b)
\end{eqnarray*}
for all $x,y \in X$, $a \in A_x$, and $b \in A_y$.  Given another right
$X^*$--module $\rmod{D} = (D,\delta,\theta)$ and a homomorphism $g\co
\rmod{B} \to \rmod{D}$ of right $X^*$--modules, we may construct another left
$X$--module $\rmod{C}$ as before, and a left $X$--module homomorphism $f\co \rmod{A} \to
\rmod{C}$ by setting $f_x = g_{x^*}$ for all $x \in X$.  This process
determines a functor $G\co \categ{RMod}^{X^*} \to \categ{RMod}_X$ which
is the inverse of $F$.
\end{proof}
\begin{corollary}
For any rack $X$, the categories $\categ{RMod}^X$ and $\categ{RMod}_{X^*}$
are equivalent
\end{corollary}
\begin{proof}
This follows immediately from the observation that $X^{**} = X$.
\end{proof}

There is a corresponding notion of a \defn{right quandle module}.  This is
a right $X$--module $\rmod{A} = (A,\phi,\psi)$ such that
$$\psi^{x,x}(a) + \phi^{x,x}(a) = a$$
for all $x \in X$ and $a \in A_x$.  We thus obtain a category
$\categ{QMod}^X$ which is equivalent to $\categ{QMod}_{X^*}$.

\section{Free $X$--modules}
For an arbitrary rack $X$, the \defn{discrete trunk} $\trunk{D}(X)$ on $X$
is the trunk with one object for each $x \in X$, and no morphisms.  A
trunk map $\rmod{S}\co \trunk{D}(X) \to \trunk{Set}$, then, determines a set
$S_x$ for each element $x \in X$.  We denote by
$\trunk{Set}^{\trunk{D}(X)}$ the category whose objects are these trunk
maps, and whose morphisms are natural transformations.

There is an obvious `forgetful' functor $U\co \categ{RMod}_X \to
\trunk{Set}^{\trunk{D}(X)}$, which, for any $X$--module $\rmod{A} =
(A,\phi,\psi)$ maps the Abelian group $A_x$ to its underlying set, and
discards the structure maps $\phi_{x,y}$ and $\psi_{y,x}$, for all $x,y
\in X$.

For any trunk map $\rmod{S}\co \trunk{D}(X) \to \trunk{Set}$, we define
the \defn{free rack module} (\defn{over} $X$) $F\rmod{S}$ to be the
module $\rmod{F} = (F,\Lambda,\Rho)$ where $F_x$ is the free abelian
group generated by symbols of the form
\begin{enumerate}
\item $\rho_{x^{\bar{w}},w}(a)$ where $w \in \As X$ and $a \in
    S_{x^{\bar{w}}}$
\item $\rho_{x^{\bar{w}},w}\lambda_{y,x^{\bar{w}\bar{y}}}(b)$
  where $w \in \As X$, $y \in X$, and $b \in S_y$
\end{enumerate}
modulo the relations
\begin{enumerate}
\item $\rho_{x^u,v}\rho_{x,u} = \rho_{x^v,u^v}\rho_{x,v} = \rho_{x,uv}$
\item $\rho_{x^y,v}\lambda_{y,x} = \lambda_{y^v,x^v}\rho_{y,v}$
\item $\lambda_{z,x^y} = \lambda_{y^z,x^z}\lambda_{z,y}
  + \rho_{x^z,y^z}\lambda_{z,x}$
\item $\rho_{x,w}(p+q) = \rho_{x,w}(p) + \rho_{x,w}(q)$
\item $\lambda_{y,x}(s+t) = \lambda_{y,x}(s) + \lambda_{y,x}(t)$
\end{enumerate}
for all $x,y,z \in X$; $u,v,w \in \As X$; $p,q \in F_x$; and $s,t \in F_y$.
The symbol $(c)$ should be interpreted as $\rho_{x,1}(c)$ for
any $c \in S_x$, with $1$ denoting the identity in $\As X$.

The structure maps are defined as follows:
\begin{eqnarray*}
\Rho_{x,y}\co F_x \to F_{x^y}; && p \mapsto \rho_{x,y}p \\
\Lambda_{y,x}\co F_y \to F_{x^y}; && s \mapsto \lambda_{y,x}s.
\end{eqnarray*}

For any two trunk maps $\rmod{S,T}\co \trunk{D}(X) \to \trunk{Set}$ and
any natural transformation $f\co \rmod{S}\to\rmod{T}$, there is a unique
$X$--map $Ff\co F\rmod{S} \to F\rmod{T}$ given by linearly extending $f$.

This functor $F\co \categ{Set}^{\trunk{D}(X)} \to \categ{RMod}_X$ is 
left adjoint to the forgetful functor $U$.

For any rack $X$, we define the \defn{rack algebra} (or \defn{wring})
$\mathbb{Z}X$ of $X$ to be the free $X$--module on the singleton trunk
map $\rmod{S}\co x \mapsto \{(*)\}$.

A typical element of $(\mathbb{Z}X)_x$ is of the form
$$\sum_{w \in \As X} n_w \rho_{x^{\bar{w}},w}(*)
  + \sum_{t \in X, v \in \As X}
  m_{t,v}\rho_{x^{\bar{v}},v}\lambda_{t,x^{\bar{v}\bar{t}}}(*)$$
where $n_w,m_{t,v} \in \mathbb{Z}$.

The composition of the structure maps in $\mathbb{Z}X$ yields a
multiplicative structure as follows:
{\small
$$\begin{array}{lclcl}
\rho_{x,v}(*) \cdot \rho_{x^{\bar{u}},u}(*) & := &
  \rho_{x,v}\rho_{x^{\bar{u}},u}(*) & = &
  \rho_{x^{\bar{u}},uv}(*) \\
\rho_{x,v}(*) \cdot \rho_{x^{\bar{u}},u}\lambda_{s,x^{\bar{u}\bar{s}}}(*)
  & := &
  \rho_{x,v}\rho_{x^{\bar{u}},u}\lambda_{s,x^{\bar{u}\bar{s}}}(*) & = &
  \rho_{x^{\bar{u}},uv}\lambda_{s,x^{\bar{u}\bar{s}}}(*) \\
\rho_{x,v}\lambda_{t,x^{\bar{t}}}(*) \cdot \rho_{t^{\bar{u}},u}(*) & := &
  \rho_{x,v}\lambda_{t,x^{\bar{t}}}\rho_{t^{\bar{u}},u}(*) & = & 
  \rho_{x^{\bar{u}},uv}\lambda_{t^{\bar{u}},x^{\bar{t}\bar{u}}}(*) \\
\rho_{x,v}\lambda_{t,x^{\bar{t}}}(*) \cdot
  \rho_{t^{\bar{u}},u}\lambda_{s,t^{\bar{u}\bar{s}}}(*) & := &
  \rho_{x,v}\lambda_{t,x^{\bar{t}}}
  \rho_{t^{\bar{u}},u}\lambda_{s,t^{\bar{u}\bar{s}}}(*)
  & = & 
  \rho_{x^{\bar{u}},uv}\lambda_{s,x^{\bar{u}\bar{s}}}(*) \\
& & & &
  - \rho_{x^{\bar{t}\bar{u}},utv}\lambda_{s,x^{\bar{t}\bar{u}\bar{s}}}(*)
\end{array}$$}
This product operation is associative and distributes over addition,
giving $\mathbb{Z}X$ a structure analogous to that of a preadditive
category or `ring with several objects'\cite{mitchell:several}.

Considering $\mathbb{Z}$ as a trivial $X$--module, we may define an
$X$--map $\epsilon\co \mathbb{Z}X \to \mathbb{Z}$ (the \defn{augmentation
map}) as follows:
$$\epsilon_x \left( \sum_{w \in \As X} n_w \rho_{x^{\bar{w}},w}(*)
  + \sum_{t \in X, v \in \As X}
  m_{t,v}\rho_{x^{\bar{v}},v}\lambda_{t,x^{\bar{v}\bar{t}}}(*) \right)
  = \sum_{w \in \As X} n_w$$
The \defn{augmentation module} of $X$ is the kernel $\augmod{X} =
\ker\epsilon$ of this map.

There are analogous constructions for quandle modules.  For any quandle
$X$, and trunk map 
$\rmod{S}\co\trunk{D}(X)\to\trunk{Set}$, the \defn{free quandle module}
(\defn{over} $X$) $F\rmod{S}$ is the free rack $X$--module $F$ on
$\rmod{S}$, modulo the relation
\begin{enumerate}
\setcounter{enumi}{5}
\item $\lambda_{x,x}(a) + \rho_{x,x}(a) = a$
\end{enumerate}
for all $x \in X$ and $a \in F_x$.  The \defn{quandle algebra} (or
\defn{wring}) of $X$, which we also denote $\mathbb{Z}X$, is the free
quandle $X$--module on the singleton trunk map $\rmod{S}\co x \mapsto
\{(*)\}$, equipped with the same multiplicative structure as the rack
algebra of $X$.

\section{Tensor products}

Given a rack $X$, let $\rmod{A} = (A,\phi,\psi)$ be a right $X$--module, and
$\rmod{B} = (B,\chi,\omega)$ be a left $X$--module.  Then the
\defn{tensor product} $\rmod{A} \otimes_X \rmod{B}$ is defined as follows.
Let $D$ be the free Abelian group with basis 
${\displaystyle\bigcup_{x \in X}}\left(A_x \times B_x\right)$.  That is,
let $D$ be generated by symbols of the form $(a,b)$ where $a \in A_x$ and
$b \in B_x$ for all $x \in X$.  Then we define $\rmod{A} \otimes_X \rmod{B}$
to be the group $D$ modulo the relations
\begin{enumerate}
\item $(a_1 + a_2, b) = (a_1,b) + (a_2,b)$
\item $(a, b_1 + b_2) = (a,b_1) + (a,b_2)$
\item $(na,b) = n(a,b) = (a,nb)$
\item $(\phi^{x,y}c,b) = (c,\chi_{x,y}b)$
\item $(\psi^{y,x}c,d) = (c,\omega_{y,x}d)$
\end{enumerate}
for all $x,y \in X$; $a,a_1,a_2 \in A_x$; $c \in A_{x^y}$; $b,b_1,b_2 \in
B_x$; $d \in B_y$; and $n \in \mathbb{Z}$.  We denote the equivalence
class of $(a,b)$ by $a \otimes b$.

An $X$--\defn{biadditive} map is an Abelian group homomorphism
$$f\co\bigcup_{x \in X}(A_x \times B_x) \to C$$
such that
\begin{enumerate}
\item $f(a_1+a_2,b) = f(a_1,b) + f(a_2,b)$
\item $f(a,b_1+b_2) = f(a,b_1) + f(a,b_2)$
\item $f(na,b) = nf(a,b) = f(a,nb)$
\item $f(\phi^{x,y}c,b) = f(c,\chi_{x,y}b)$
\item $f(\psi^{y,x}c,d) = f(c,\omega_{y,x}d)$
\end{enumerate}
for all $x,y \in X$; $a,a_1,a_2 \in A_x$; $c \in A_{x^y}$; $b,b_1,b_2 \in
B_x$; $d \in B_y$; and $n \in \mathbb{Z}$.  The tensor product $\rmod{A}
\otimes_X \rmod{B}$, then, has the universal property that, for any
Abelian group $C$ and $X$--biadditive map $f\co{\displaystyle\bigcup_{x
\in X}}\left(A_x \times B_x\right) \to C$, there is a unique Abelian
group homomorphism $g\co\rmod{A} \otimes_X \rmod{B}\to C$ making the
diagram
$$\Vtriangle[{\displaystyle\bigcup_{x \in X}}\left(A_x \times B_x\right)
            `\rmod{A}\otimes_X\rmod{B}`C;h`f`g]$$
commute.  The tensor product is unique up to isomorphism by the usual
universality argument.

We may give the set $\Hom_{\categ{Ab}}(\rmod{B},C)$ a canonical right
$X$--module structure.  Let $H_x = \Hom_{\categ{Ab}}(B_x,C)$, which has an
obvious Abelian group structure defined by $(f_x + g_x)(b) := f_x(b) +
g_x(b)$ for all $f_x, g_x\co B_x \to C$.  Now define
structure maps $\eta^{x,y}\co H_{x^y} \to H_x$ and $\zeta^{y,x}\co H_{x^y}
\to H_y$ by
\begin{eqnarray*}
\eta^{x,y}(f) & = & f\chi_{x,y}\co B_x \to C \\
\zeta^{y,x}(f) & = & f\omega_{y,x}\co B_y \to C
\end{eqnarray*}
for all $x,y \in X$, and $f\co B_{x^y} \to C$.  Then $\rmod{H} =
(H,\eta,\zeta)$ is a right $X$--module, and this construction gives a
functor $H_{\rmod{B}}(-) = \Hom_{\categ{Ab}}(\rmod{B},-)\co
\categ{Ab}\to\categ{RMod}^X$.

\begin{proposition}
For any right $X$--module $\rmod{A} = (A,\phi,\psi)$, left $X$--module
$\rmod{B} = (B,\chi,\omega)$, and abelian group $C$, there is a natural
isomorphism
$$\Hom_{\categ{Ab}}(\rmod{A} \otimes_X \rmod{B},C)
  \iso \Hom_{\categ{RMod}^X}(\rmod{A},\Hom_{\categ{Ab}}(\rmod{B},C)).$$
That is, the functor $-\otimes_X\rmod{B}\co\categ{RMod}^X\to\categ{Ab}$ is
left adjoint to the functor
$H_{\rmod{B}} = \Hom_{\categ{Ab}}(\rmod{B},-)\co\categ{Ab}\to\categ{RMod}^X$.
\end{proposition}

\begin{proof}
An element $f$ of
$\Hom_{\categ{RMod}^X}(\rmod{A},\Hom_{\categ{Ab}}(\rmod{B},C))$ assigns,
to each $x \in X$ and $a \in A_x$, an Abelian group homomorphism
$f_x(a)\co B_x \to C$, in a natural way.  That is,
$$
\begin{array}{lclcl}
f_x(a_1+a_2)(b) & = & f_x(a_1)(b) + f_x(a_2)(b) & & \\
f_x(a)(b_1+b_2) & = & f_x(a)(b_1) + f_x(a)(b_2) & & \\
f_x(na)(b)      & = & nf_x(a)(b) & = & f_x(a)(nb) \\
f_x(\phi^{x,y}(c)(b) & = & \left(\eta^{y,x}f_{x^y}(c)\right)(b) & = &
  f_{x^y}(c)(\chi_{x,y}b) \\
f_y(\psi^{x,y}(c)(d) & = & \left(\zeta^{y,x}f_{x^y}(c)\right)(d) & = &
  f_{x^y}(c)(\omega_{y,x}d)
\end{array}$$
for all $x,y \in X$; $a,a_1,a_2 \in A_x$; $c \in A_{x^y}$; $b,b_1,b_2 \in
B_x$; $d \in B_y$; and $n \in \mathbb{Z}$.  The first three identities follow
from the fact that each $f_x$ is an abelian group homomorphism, and the
remaining two from the fact that $f$ is a homomorphism of right
$X$--modules (and hence a natural transformation of contravariant trunk
maps $\trunk{T}(X)\to\trunk{Ab}$).

The required natural isomorphism is
$$\tau_{\rmod{A},C}\co
\Hom_{\categ{Ab}}(\rmod{A}\otimes_X\rmod{B},C) \to
\Hom_{\categ{RMod}^X}(\rmod{A},\Hom_{\categ{Ab}}(\rmod{B},C))$$
defined by
$$f(a \otimes b) \mapsto f_x(a)(b)$$
for all $x \in X, a \in A_x$, and $b \in B_x$.
\end{proof}

Although the tensor product of two $X$--modules is an Abelian group, in
certain circumstances it may itself be regarded as an $X$--module in a
canonical way.
\begin{proposition}
For any left $X$--module $\rmod{A} = (A,\phi,\psi)$, the tensor product
$\mathbb{Z}X\otimes_X\rmod{A}$ has a canonical left $X$--module structure
such that $\mathbb{Z}X\otimes_X\rmod{A}\iso\rmod{A}$.
\end{proposition}

\begin{proof}
The tensor product $\mathbb{Z}X \otimes_X \rmod{A}$ is the free Abelian
group generated by symbols of the form $r \in (\mathbb{Z}X)_x$ (where
$\mathbb{Z}X$ is considered as a right $X$--module) and $a \in A_x$, for
all $x \in X$, such that
\begin{enumerate}
\item $\rho^{x,y}(*) \otimes a = (*) \otimes \phi_{x,y}(a)$
\item $\lambda^{y,x}(*) \otimes b = (*) \otimes \psi_{y,x}(b)$
\end{enumerate}
where $a \in A_x$ and $b \in A_y$.  For each $x \in X$, let $B_x$ be the
free Abelian group generated by symbols of the form $(*)\otimes a$ where
$a \in A_x$.  We may then define homomorphisms
\begin{eqnarray*}
\chi_{x,y}\co B_x \to B_{x^y} & ; &
  (*) \otimes a \mapsto \rho^{x,y}(*) \otimes a = (*) \otimes \phi_{x,y}(a) \\
\omega_{y,x}\co B_y \to B_{x^y} & ; &
  (*) \otimes b \mapsto \lambda^{x,y}(*) \otimes b = (*) \otimes \psi_{y,x}(b)
\end{eqnarray*}
which satisfy the criteria for the structure maps of a left $X$--module.
This $X$--module $\rmod{B} = (B,\chi,\omega)$ is isomorphic to $\rmod{A}$.
\end{proof}

\section{Homology and cohomology}
Given a rack $X$ and an $X$--module $\rmod{A} = (A,\phi,\psi)$, the group
$\Ext(X,\rmod{A})$ (defined in \cite{jackson:extrq}) is the group of
equivalence classes of (factor sets corresponding to) extensions of $X$ by
$\rmod{A}$.  This may be defined to be the quotient
$Z(X,\rmod{A})/B(X,\rmod{A})$, where $Z(X,\rmod{A})$ consists of factor
sets $\sigma$ satisfying the condition
$$\sigma_{x^y,z} + \phi_{x^y,z}(\sigma_{x,y}) =
  \psi_{y^z,x^z}(\sigma_{y,z}) + \sigma_{x^z,y^z} +
  \phi_{x^z,y^z}(\sigma_{x,z})$$
and $B(X,\rmod{A})$ is the subgroup of $Z(X,\rmod{A})$ consisting of
factor sets of the form
$$\sigma_{x,y} = \psi_{y,x}(\upsilon_y) - \upsilon_{x^y} +
  \phi_{x,y}(\upsilon_x)$$
for all $x,y,z \in X$, and some set $\upsilon = \{\upsilon_x \in A_x : x
\in X\}$
In the case where $\rmod{A}$ is a trivial $X$--module, this reduces to the
definition of the second cohomology $H^2(X;A)$ of (the rack space of) $X$
with coefficients in the Abelian group $A$.  We wish to formalise this
connection with rack cohomology, and devise suitable generalisations of
the higher homology and cohomology groups of $X$ to the case where the
coefficient object is an arbitrary $X$--module, rather than just an
Abelian group.

We begin by considering the elements of $B(X,\rmod{A})$, expecting them to
be the image, under some suitable coboundary operator, of `functions'
$f\co X \to \rmod{A}$.  This concept is not yet well-defined, as it is not
immediately clear what is meant by a `function' from a rack (which is a
set with some additional structure imposed on it) to a rack module (which
is a trunk map).

Initially, then, we define a 1--\defn{coboundary} to be a family $\upsilon
= \{\upsilon_x \in A_x : x \in X\}$ of group elements such that
$$\upsilon_{x^y} = \psi_{y,x}(\upsilon_y) + \phi_{x,y}(\upsilon_x).$$
In order to reformulate this notion in a more useful manner, we must
find some way of describing the rack $X$ as a trunk map $\trunk{D}(X) \to
\trunk{Set}$.

Let $\rmod{S}_1$ denote the trunk map $\trunk{D}(X) \to \trunk{Set}$ where
$(\rmod{S}_1)_x = \{x\}$ for all $x \in X$.  Then a 1--coboundary
$\upsilon$ is a natural transformation
$\upsilon\co\rmod{S}_1\to U\rmod{A}$, where $U$ denotes the forgetful
functor $\categ{RMod}_X\to\trunk{Set}^{\trunk{D}(X)}$.  This set
$\Hom_{\trunk{Set}^{\trunk{D}(X)}}(\rmod{S}_1,U\rmod{A})$ has an
Abelian group structure defined by setting $(\upsilon + \omega)_x(x) =
\upsilon_x(x) + \omega_x(x)$ for all $x \in X$.

Similarly, define $\rmod{S}_2\co \trunk{D}(X) \to \trunk{Set}$ by
$(\rmod{S}_2)_x = \{(p,q) \in X \times X : p^q = x \}$.  A factor set
$\sigma$ may be regarded as a natural transformation
$\sigma\co\rmod{S}_2\to U\rmod{A}$.  The set
$\Hom_{\trunk{Set}^{\trunk{D}(X)}}(\rmod{S}_2,U\rmod{A})$ also has an
obvious Abelian group structure, defined by setting $(\sigma +
\tau)_{p^q}(p,q) = \sigma_{p^q}(p,q) + \tau_{p^q}(p,q)$ for all $p,q \in
X$.

If we now define a map $d^2\co
\Hom_{\trunk{Set}^{\trunk{D}(X)}}(\rmod{S}_1,U\rmod{A}) \to
\Hom_{\trunk{Set}^{\trunk{D}(X)}}(\rmod{S}_2,U\rmod{A})$ by
$$(d^2f)_{x^y}(x,y) = \psi_{x,y}f_y(y) - f_{x^y}(x^y) + \phi_{x,y}f_x(x)$$
then $B(X,\rmod{A})$ may be seen to be exactly $\im d^2$.  This map $d^2$
is an Abelian group homomorphism.

In general, define the trunk map $\rmod{S}_n\co \trunk{D}(X) \to
\trunk{Set}$ by
$$(\rmod{S}_n)_x = \{ (x_1,\ldots,x_n) \in X^n : x_1^{x_2 \ldots x_n} = x \}$$
for $n>0$, and
$$(\rmod{S}_0)_x = \{(*)\}$$
for all $x \in X$.

Let $d^3\co \Hom_{\trunk{Set}^{\trunk{D}(X)}}(\rmod{S}_2,U\rmod{A}) \to
\Hom_{\trunk{Set}^{\trunk{D}(X)}}(\rmod{S}_3,U\rmod{A})$ such that
\begin{multline*}
(d^3f)_{x^{yz}}(x,y,z) = \psi_{y^z,x^z}f_{y^z}(y,z) \\
  + f_{x^{zy^z}}(x^z,y^z) - f_{x^{yz}}(x^y,z) \\
  + \phi_{x^z,y^z}f_{x^z}(x,z) - \phi_{x^y,z}f_{x^y}(x,y).
\end{multline*}
Then $Z(X,\rmod{A}) = \ker d^3$.  A routine calculation confirms that $\im
d^2 \leqslant \ker d^3$, and so we now have a fragment
$$\Hom_{\trunk{Set}^{\trunk{D}(X)}}(\rmod{S}_1,U\rmod{A})
  \stackrel{d^2}{\longrightarrow}
  \Hom_{\trunk{Set}^{\trunk{D}(X)}}(\rmod{S}_2,U\rmod{A})
  \stackrel{d^3}{\longrightarrow}
  \Hom_{\trunk{Set}^{\trunk{D}(X)}}(\rmod{S}_3,U\rmod{A})$$
of a cochain complex of Abelian groups.  This is equivalent to
$$\Hom_X(\rmod{F}_1,\rmod{A})
  \stackrel{d^2}{\longrightarrow}
  \Hom_X(\rmod{F}_2,\rmod{A})
  \stackrel{d^3}{\longrightarrow}
  \Hom_X(\rmod{F}_3,\rmod{A})$$
where $\rmod{F}_n$ is the free $X$--module with basis $\rmod{S}_n$.  We
thus have a description of the second cohomology of $X$ in terms of the
application of the contravariant functor $\Hom_X(-,\rmod{A})$ to part of a
complex of free $X$--modules.

Seeking a similar perspective for the first cohomology $H^1(X;\rmod{A})$,
we define a \defn{derivation} $f\co X \to \rmod{A}$ to be a natural
transformation $f\co\rmod{S}_1 \to U\rmod{A}$ such that
$$f_{x^y}(x^y) = \psi_{y,x}f_y(y) + \phi_{x,y}f_x(x);$$
that is, an element of $\ker d^2$.  We denote the group of such maps by
$\Der(X,\rmod{A})$.  If $\rmod{A}$ is a trivial $X$--module, then
$\Der(X,\rmod{A}) = \Hom_X(\rmod{F}_1,\rmod{A})$.

If $z$ is an arbitrary fixed element of $X$, then a $z$--\defn{principal
derivation} is a natural transformation $f\co\rmod{S}_1 \to U\rmod{A}$
such that $f_x(x) = \psi_{z,x^{\bar z}}(a)$ for each $x \in X$ and some
fixed element $a \in A_z$.  As the terminology suggests, such a map is
itself a derivation.  Furthermore, the set $\PDer(X,\rmod{A})$ of all
$z$--principal derivations (of $X$ into $\rmod{A}$) has an Abelian group
structure, and may be regarded as the image of the $X$--map $d^1_z\co
\Hom_{\trunk{Set}^{\trunk{D}(X)}}(\rmod{S}_0,U\rmod{A}) \to
\Hom_{\trunk{Set}^{\trunk{D}(X)}}(\rmod{S}_1,U\rmod{A})$ given by
$$(d^1_zf)_x(x) = \psi_{z,x^{\bar z}}(*).$$
We may define $H^1(X;\rmod{A})_z$, the first cohomology group of $X$ (with
respect to $z$) with coefficients in $\rmod{A}$, to be the quotient
$\Der(X,\rmod{A})/\PDer(X,\rmod{A})$.

We have thus extended our cochain complex fragment by one dimension:
$$\Hom_X(\rmod{F}_0,\rmod{A})
  \stackrel{d^1_z}{\longrightarrow}
  \Hom_X(\rmod{F}_1,\rmod{A})
  \stackrel{d^2}{\longrightarrow}
  \Hom_X(\rmod{F}_2,\rmod{A})
  \stackrel{d^3}{\longrightarrow}
  \Hom_X(\rmod{F}_3,\rmod{A}).$$
In many cases (when $\rmod{A}$ is a homogeneous trivial, dihedral, or
Alexander module, for example) the first cohomology is independent of
the choice of fixed element of $X$, and we may omit the $z$ subscript.
In particular, when $\rmod{A}$ is a trivial module, the group
$\PDer(X,\rmod{A})$ is itself trivial.

With this discussion in mind, we now define the ($z$--)\defn{standard
complex} of a rack $X$ to be
$$\mathbf{F}_z = \cdots \stackrel{d_2}{\longrightarrow} \rmod{F}_1
  \stackrel{d_1^z}{\longrightarrow} \rmod{F}_0
  \stackrel{\varepsilon}{\longrightarrow} \mathbb{Z}
  \longrightarrow 0$$
where $\rmod{F}_n$ is the free $X$--module with basis $\rmod{S}_n$,
and the boundary maps are given by
$$d_n = \sum_{i=1}^n (-1)^{i+1}d_n^i$$
where
\begin{multline*}
(d_n^i)_{x_1^{x_2 \ldots x_n}}(x_1,\ldots,x_n)  = \\
  \rho_{x_1^{x_2 \ldots \widehat{x_{i+1}}\ldots x_n},
          x_{i+1}^{x_{i+2} \ldots x_n}}
    (x_1,\ldots,\widehat{x_{i+1}},\ldots,x_n) \\
  - (x_1^{x_{i+1}},\ldots,x_i^{x_{i+1}},x_{i+2},\ldots,x_n)
\end{multline*}
for $1 \leqslant i \leqslant n-1$,
\begin{equation*}
(d_n^n)_{x_1^{x_2 \ldots x_n}}(x_1,\ldots,x_n) =
  (-1)^{n+1}\lambda_{x_2^{x_3 \ldots x_n}, x_1^{x_3 \ldots x_n}}
(x_2,\ldots,x_n)
\end{equation*}
for $n > 1$, and
\begin{equation*}
(d_1^z)_x\co (x) \mapsto \lambda_{z,x^{\bar{z}}}(*)
\end{equation*}
where $\widehat{.}$ denotes elision of the marked symbol, and where $z$ is
an
arbitrarily chosen, fixed element of $X$.  In particular, both
$\rmod{F}_1$ and $\rmod{F}_0$ are isomorphic to the rack algebra
$\mathbb{Z} X$.
We set the map $d_0\co \rmod{F}_0 \to \mathbb{Z}$ to be the augmentation map
$\epsilon\co \mathbb{Z} X \to \mathbb{Z}$.
\begin{lemma}
If $1 \leqslant i < j < n$ then
$$d_{n-1}^i d_n^j = d_{n-1}^{j-1} d_n^i.$$
If $1 \leqslant i < n-1$ then
$$d_{n-1}^i d_n^n = - d_{n-1}^{n-1}d_n^{i+1}.$$
And
$$d_{n-1}^{n-1}d_n^n = (-1)^n d_{n-1}^{n-1}d_n^1$$
\end{lemma}
\begin{proof}
If $1 \leqslant i < j < n$,
\begin{multline*}
(d^i_{n-1})_{x_1^{x_2 \ldots x_n}}(d^j_n)_{x_1^{x_2 \ldots x_n}}
  (x_1,\ldots,x_n) = \\
\rho_{x_1^{x_2 \ldots \widehat{x_{j+1}} \ldots x_n},
      x_{j+1}^{x_{j+2} \ldots x_n}}
  \rho_{x_1^{x_2 \ldots \widehat{x_{i+1}} \ldots \widehat{x_{j+1}} \ldots
x_n},
        x_{i+1}^{x_{i+2} \ldots x_n}} \\
  (x_1,\ldots,\widehat{x_{i+1}},\ldots,\widehat{x_{j+1}},\ldots,x_n) \\
- \rho_{x_1^{x_2 \ldots \widehat{x_{j+1}} \ldots x_n},
        x_{j+1}^{x_{j+2} \ldots x_n}}
  (x_1^{x_{i+1}},\ldots,x_i^{x_{i+1}},x_{i+2},\ldots,x_n) \\
- \rho_{x_1^{x_2 \ldots \widehat{x_{i+1}} \ldots x_n},
      x_{i+1}^{x_{i+2} \ldots x_n}}
  (x_1^{x_{j+1}},\ldots,x_j^{x_{j+1}},x_{j+2},\ldots,x_n) \\
+ (x_1^{x_{i+1}x_{j+1}},\ldots,x_i^{x_{i+1}x_{j+1}},x_{i+2}^{x_{j+1}},
   \ldots,x_j^{x_{j+1}},x_{j+2},\ldots,x_n) \\
= (d^{j-1}_{n-1})_{x_1^{x_2 \ldots x_n}}(d^i_n)_{x_1^{x_2 \ldots x_n}}
  (x_1,\ldots,x_n).
\end{multline*}
If $1 \leqslant i < n-1$,
\begin{multline*}
(d^i_{n-1})_{x_1^{x_2 \ldots x_n}}(d_n^n)_{x_1^{x_2 \ldots x_n}}
  (x_1,\ldots,x_n) = \\
= (-1)^{n+1}\lambda_{x_2^{x_3 \ldots x_n},x_1^{x_3 \ldots x_n}}
  \rho_{x_2^{x_3 \ldots \widehat{x_{i+1}} \ldots x_n},
        x_{i+1}^{x_{i+2} \ldots x_n}}
  (x_2,\ldots,\widehat{x_{i+1}},\ldots,x_n) \\
- (-1)^{n+1}\lambda_{x_2^{x_3 \ldots x_n},x_1^{x_3 \ldots x_n}}
  (x_2^{x_{i+1}},\ldots,x_i^{x_{i+1}},x_{i+2},\ldots,x_n) \\
= -(d^{n-1}_{n-1})_{x_1^{x_2 \ldots x_n}}
  (d^{i+1}_n)_{x_1^{x_2 \ldots x_n}}(x_1,\ldots,x_n).
\end{multline*}
Finally,
\begin{multline*}
(d_{n-1}^{n-1})_{x_1^{x_2 \ldots x_n}}
  (d_n^n)_{x_1^{x_2 \ldots x_n}}(x_1,\ldots,x_n) = \\
(-1)^{2n+1}\lambda_{x_2^{x_3 \ldots x_n},x_1^{x_3 \ldots x_n}}
  \lambda_{x_3^{x_4 \ldots x_n},x_2^{x_4 \ldots x_n}}(x_3,\ldots,x_n) = \\
(-1)^{2n}\rho_{x_1^{x_3 \ldots x_n},x_2^{x_3 \ldots x_n}}
  \lambda_{x_3^{x_4 \ldots x_n},x_1^{x_4 \ldots x_n}}(x_3,\ldots,x_n) \\
  - (-1)^{2n}\lambda_{x_3^{x_4 \ldots x_n},x_1^{x_2 x_4 \ldots x_n}}
  (x_3,\ldots,x_n) = \\
(-1)^n (d^{n-1}_{n-1})_{x_1^{x_2 \ldots x_n}}(d_n^1)_{x_1^{x_2 \ldots
x_n}}
  (x_1,\ldots,x_n).
\end{multline*}
\end{proof}

\begin{theorem}
The standard complex is indeed a chain complex of $X$--modules.
\end{theorem}
\begin{proof}
Using the above lemma, we find that
\begin{multline*}
d_{n-1}d_n = \sum_{i=1}^{n-1}\sum_{j=1}^n (-1)^{i+j} d_{n-1}^i d_n^j = \\
\sum_{i=1}^{n-2}\sum_{j=1}^i (-1)^{i+j} d_{n-1}^i d_n^j
+ \sum_{i=1}^{n-2}\sum_{j=i+1}^{n-1} (-1)^{i+j} d_{n-1}^{j-1} d_n^i \\
+ \sum_{i=1}^{n-2} (-1)^{i+n} d_{n-1}^i d_n^n
+ \sum_{i=1}^{n-2} (-1)^{i+n} d_{n-1}^{n-1} d_n^{i+1} \\
+ (-1)^n d_{n-1}^{n-1} d_n^n
+ (-1)^{2n-1} d_{n-1}^{n-1} d_n^1 = 0
\end{multline*}
as asserted.
\end{proof}

We may now use this complex to define the homology and cohomology groups of
a rack $X$, with coefficients in an arbitrary $X$--module $\rmod{A}$:
\begin{eqnarray*}
H_n(X;\rmod{A})_z & = & H_n(\mathbf{F}_z \otimes_X \rmod{A}) \\
H^n(X;\rmod{A})_z & = & H^n(\Hom_X(\mathbf{F}_z,\rmod{A}))
\end{eqnarray*}
for $n\geqslant0$.

In the case where the coefficient module is a trivial homogeneous
$X$--module (equivalently, an Abelian group), this homology and
cohomology theory is equivalent to the (topological) homology and
cohomology of the rack space $BX$, as introduced by Fenn, Rourke and
Sanderson \cite{frs:trunks}.

If the coefficient module $\rmod{A} = (A,\phi,\psi)$ is homogeneous, then
we recover the rack homology and cohomology theories described by
Andruskiewitsch and Gra\~na \cite{ag:pointed}, and if, in addition, the
$\psi$--maps are zero (giving $\rmod{A}$ the structure of an $\As
X$--module) then we recover the cohomology theory studied by Etingof and
Gra\~na \cite{eg:orc}.

We now investigate the specialisation to the subcategory
$\categ{Quandle}$.  The construction of the standard complex
$\mathbf{F}^{\scriptscriptstyle Q}_z$ for quandle homology and cohomology
is very similar to that for rack homology and cohomology, with $\rmod{F}_n
= F\rmod{S}_n$ where
$$(\rmod{S}_n)_x = \{(x_1,\ldots,x_n) \in X^n : x_1^{x_2 \ldots x_n} = x,
\text{and}~ x_i \not= x_{i+1} ~ \text{for} ~ 1 \leqslant i \leqslant n-1 \}$$
and the same boundary maps $d_n\co \rmod{F}_n \to \rmod{F}_{n-1}$.  This
complex may be used to define quandle homology and cohomology groups
\begin{eqnarray*}
H^{\scriptscriptstyle Q}_n(X;\rmod{A})_z & = &
  H_n(\mathbf{F}^{\scriptscriptstyle Q}_z \otimes_X \rmod{A}) \\
H_{\scriptscriptstyle Q}^n(X;\rmod{A})_z & = &
  H^n(\Hom_X(\mathbf{F}^{\scriptscriptstyle Q}_z,\rmod{A}))
\end{eqnarray*}
for $n\geqslant0$.

If the coefficient module is a trivial homogeneous quandle $X$--module
(equivalently, an Abelian group), this theory is equivalent to
the one introduced by Carter, Jelsovsky, Kamada, Langford and
Saito\cite{cjkls:qcoh-state}, and if the coefficient module is a
homogeneous Alexander module \cite[Example 2.4]{jackson:extrq},
we recover Carter, Elhamdadi and Saito's twisted quandle homology
theory\cite{ces:twisted}.

\end{document}